\documentclass[12pt]{article}
\usepackage{amsmath, amsthm, amsfonts, amscd, amssymb, latexsym}
\usepackage{color}
\usepackage{enumitem}
\setlist[enumerate,1]{%
label={\normalfont(\arabic*)},ref=\arabic*%
}

\theoremstyle{plain}
\newtheorem{thm}{Theorem}[section]

\newtheorem{prop}[thm]{Proposition}

\theoremstyle{definition}
\newtheorem{defn}[thm]{Definition}
\newtheorem{ex}[thm]{Example}

\numberwithin{equation}{section}
\begin{document}
\title{\textbf{Biharmonic Hermitian vector bundles over compact 
K\"{a}hler Einstein manifolds}}
\author{
Hajime Urakawa
\\
{\small Tohoku University, Graduate School of Information Sciences, 
Division of Mathematics}
\\
{\small Aoba 6-3-09, Sendai 980-8579, Japan, urakawa@math.is.tohoku.ac.jp }}
\maketitle
\begin{abstract} In this paper, we show that, for every Hermitian vector bundle 
$\pi:\,(E,g)\rightarrow (M,h)$ over 
a compact K\"{a}hler Einstein manifold $(M,h)$,  
if the projection $\pi$ is biharmonic, then it is harmonic. 
\end{abstract}
\section{Introduction.} 
Research of harmonic maps, which are critical points of the energy functional, 
is one of the central problems in differential geometry including  minimal submanifolds.
The Euler-Lagrange equation is given by the vanishing of the tension field. 
In 1983, Eells and Lemaire (\cite{EL}) proposed 
to study biharmonic maps, which are critical points of the bienergy functional, 
by definition,
half of the integral of square of the norm of tension field $\tau(\varphi)$ 
for a smooth map $\varphi$ of a Riemannian manifold $(M,g)$
into another Riemannian manifold $(N,h)$. 
After a work by G.Y. Jiang \cite{J},
several geometers have studied biharmonic maps 
(see \cite{CMP}, \cite{IIU1}, \cite{IIU2}, \cite{II}, 
\cite{LO}, \cite{MO}, \cite{OT2}, \cite{S}, etc.). 
Note that a harmonic maps is always biharmonic. 
One of central problems is to ask whether the converse is true. 
{\em B.-Y. Chen's conjecture} is 
to ask whether every biharmonic submanifold of the Euclidean space 
${\mathbb R}^n$ must be harmonic, i.e., minimal (\cite{C}).  
There are many works supporting this conjecture (\cite{D}, \cite{HV}, \cite{KU}, \cite{AM}).
However, B.-Y. Chen's conjecture is still open.
R. Caddeo, S. Montaldo, P. Piu (\cite{CMP}) and C. Oniciuc (\cite{Oniciuc})
raised the {\em generalized B.-Y. Chen's conjecture} to ask
whether each biharmonic submanifold in a Riemannian manifold $(N,h)$
of non-positive sectional curvature must be harmonic (minimal).
For the generalized Chen's conjecture,
Ou and Tang gave (\cite{OT1}, \cite{OT2}) a counter example in some Riemannian manifold
of negative sectional curvature. 
But, it is also known (cf. \cite{NU1}, \cite{NU2}, \cite{NUG})
that every biharmonic map of a complete Riemannian manifold
into another Riemannian manifold of non-positive sectional curvature
with finite energy and finite bienergy must be harmonic. 
For the target Riemannian manifold $(N,h)$ of non-negative sectional curvature, 
theories of biharmonic maps and biharmonic immersions seems to be quite different
from the case $(N,h)$ of non-positive sectional curvature. 
There exit biharmonic submanifolds which is not harmonic 
in the unit sphere.
S. Ohno, T. Sakai and myself \cite{OSU1}, \cite{OSU2} 
determined (1) all the biharmonic hypersurfaces in irreducible symmetric spaces of compact type
which are regular orbits of commutative Hermann actions of cohomogeneity one, and gave 
 (2) a complete table of all the proper biharmonic singular orbits
of commutative Hermann actions of cohomogeneity two,
and (3) a complete list of all the proper biharmonic regular orbits of
$(K_2 \times K_1)$-actions of cohomogeneity one on $G$ for every  commutative compact symmetric triad
$(G, K_1, K_2)$.
We note that recently Inoguchi and Sasahara (\cite{IS}) also investigated biharmonic homogeneous hypersurfaces
in compact symmetric spaces,
and Ohno studied biharmonic orbits of isotropy representations of symmetric spaces
in the sphere (cf. \cite{Oh1}, \cite{Oh2}).
\vskip0.6cm\par
In  this paper, 
we treat with an Hermitian vector bundle $(E,g)\rightarrow (M,h)$ 
over a compact Riemannian manifold $(M,h)$. 
We assume $(M,h)$ is a compact K\"{a}hler Einstein Riemannian manifold, 
that is, 
 the Ricci transform ${\rm Ric}^h$ of 
 the K\"{a}hler metric $h$ on $M$ satisfies 
${\rm Ric}^h=c\,\,{\rm Id}$, for some constant $c$. 
Then, we show the following: 
\begin{thm}
Let $\pi:\,\,(E,g)\rightarrow (M,h)$ be an Hermitian vector bundle 
over a compact K\"{a}hler Einstein Riemannian manifold 
$(M,h)$. If $\pi$ is biharmonic, then it is harmonic. 
\end{thm}
\vskip0.6cm\par 
Theorem 1.1 shows the sharp contrasts on the biharmonicities 
between the case of vector bundles and 
the one of the principle $G$-bundles. 
Indeed, we treated with the biharmonicity of 
the projection of the 
principal $G$-bundle over a Riemannian manifold $(M,h)$ 
with negative definite Ricci tensor field (cf. Theorem 2.3 in \cite{U4}).  
We also gave an example of 
the projection of the 
principal $G$-bundle over a Riemannian manifold $(M,h)$ 
which is biharmonic but not harmonic 
(cf. Theorem 5 in \cite{U5}). 
\vskip0.3cm\par
\section{Preliminaries.}
In this section, we prepare necessary materials for the first and second variational formulas for the bienergy functional and biharmonic maps. 
Let us recall the definition of a harmonic map $\varphi:\,(M,g)\rightarrow (N,h)$, of a compact Riemannian manifold $(M,g)$ into a Riemannian manifold $(N,h)$, 
which is an extremal 
of the {\em energy functional} defined by 
$$
E(\varphi)=\int_Me(\varphi)\,v_g, 
$$
where $e(\varphi):=\frac12\vert d\varphi\vert^2$ is called the energy density 
of $\varphi$.  
That is, for any variation $\{\varphi_t\}$ of $\varphi$ with 
$\varphi_0=\varphi$, 
\begin{equation}
\frac{d}{dt}\bigg\vert_{t=0}E(\varphi_t)=-\int_Mh(\tau(\varphi),V)v_g=0,
\end{equation}
where $V\in \Gamma(\varphi^{-1}TN)$ is the variation vector field along $\varphi$ which is given by 
$V(x)=\frac{d}{dt}\vert_{t=0}\varphi_t(x)\in T_{\varphi(x)}N$, 
$(x\in M)$, 
and  the {\em tension field} is given by 
$\tau(\varphi)
=\sum_{i=1}^mB(\varphi)(e_i,e_i)\in \Gamma(\varphi^{-1}TN)$, 
where 
$\{e_i\}_{i=1}^m$ is a locally defined frame field on $(M,g)$, 
and $B(\varphi)$ is the second fundamental form of $\varphi$ 
defined by 
\begin{align}
B(\varphi)(X,Y)&=(\widetilde{\nabla}d\varphi)(X,Y)\nonumber\\
&=(\widetilde{\nabla}_Xd\varphi)(Y)\nonumber\\
&=\overline{\nabla}_X(d\varphi(Y))-d\varphi(\nabla_XY)\nonumber\\
&=\nabla^h_{d\varphi(X)}d\varphi(Y)-d\varphi(\nabla_XY),
\end{align}
for all vector fields $X, Y\in {\frak X}(M)$. 
Furthermore, 
$\nabla$, and
$\nabla^h$, 
 are the Levi-Civita connections on $TM$, $TN$  of $(M,g)$, $(N,h)$, respectively, and 
$\overline{\nabla}$, and $\widetilde{\nabla}$ are the induced ones on $\varphi^{-1}TN$, and $T^{\ast}M\otimes \varphi^{-1}TN$, respectively. By (2.1), $\varphi$ is harmonic if and only if $\tau(\varphi)=0$. 
\par
The second variation formula is given as follows. Assume that 
$\varphi$ is harmonic. 
Then, 
\begin{equation}
\frac{d^2}{dt^2}\bigg\vert_{t=0}E(\varphi_t)
=\int_Mh(J(V),V)v_g, 
\end{equation}
where 
$J$ is an elliptic differential operator, called 
{\em Jacobi operator}  acting on 
$\Gamma(\varphi^{-1}TN)$ given by 
\begin{equation}
J(V)=\overline{\Delta}V-{\mathcal R}(V),
\end{equation}
where 
$\overline{\Delta}V=\overline{\nabla}^{\ast}\overline{\nabla}V
=-\sum_{i=1}^m\{
\overline{\nabla}_{e_i}\overline{\nabla}_{e_i}V-\overline{\nabla}_{\nabla_{e_i}e_i}V
\}$ 
is the {\em rough Laplacian} and 
${\mathcal R}$ is a linear operator on $\Gamma(\varphi^{-1}TN)$
given by 
${\mathcal R}(V)=
\sum_{i=1}^mR^h(V,d\varphi(e_i))d\varphi(e_i)$,
and $R^h$ is the curvature tensor of $(N,h)$ given by 
$R^h(U,V)=\nabla^h{}_U\nabla^h{}_V-\nabla^h{}_V\nabla^h{}_U-\nabla^h{}_{[U,V]}$ for $U,\,V\in {\frak X}(N)$.   
\par
J. Eells and L. Lemaire \cite{EL} 
proposed polyharmonic ($k$-harmonic) maps and 
Jiang \cite{J} studied the first and second variation formulas of biharmonic maps. Let us consider the {\em bienergy functional} 
defined by 
\begin{equation}
E_2(\varphi)=\frac12\int_M\vert\tau(\varphi)\vert ^2v_g, 
\end{equation}
where 
$\vert V\vert^2=h(V,V)$, $V\in \Gamma(\varphi^{-1}TN)$.  
\par
Then, the first variation formula of the bienergy functional 
is given as follows. 
\begin{thm}
\quad $($the first variation formula$)$ 
\begin{equation}
\frac{d}{dt}\bigg\vert_{t=0}E_2(\varphi_t)
=-\int_Mh(\tau_2(\varphi),V)v_g.
\end{equation}
Here, 
\begin{equation}
\tau_2(\varphi)
:=J(\tau(\varphi))=\overline{\Delta}\tau(\varphi)-{\mathcal R}(\tau(\varphi)),
\end{equation}
which is called the {\em bitension field} of $\varphi$, and 
$J$ is given in $(2.4)$.  
\end{thm}
\begin{defn}
A smooth map $\varphi$ of $M$ into $N$ is said to be 
{\em biharmonic} if 
$\tau_2(\varphi)=0$. 
\end{defn} 
\vskip0.6cm\par
\section{Proof of Theorem 1.1.}
\vskip0.6cm\par
To prove Theorem 1.1, we need the following: 
\begin{prop}Let $\pi:\,\,(E,g)\rightarrow (M,g)$ be an Hermitian vector bundle 
over a compact K\"{a}hler Einstein manifold 
$(M,h)$. 
Assume that $\pi$ is biharmonic. Then the following hold: 
\par
$(1)$ The tension field 
$\tau(\pi)$ satisfies that 
\begin{equation}
\overline{\nabla}_X\tau(\pi)=0\qquad\qquad\qquad\qquad 
 (\forall \,\,X\in {\frak X}(M)).
\end{equation}
\par
$(2)$ The pointwise inner product 
$\langle \tau(\pi),\tau(\pi)\rangle=\vert\,\tau(\pi)\,\vert^2$ is constant on $(M,g)$, say $d\geq 0$. 
 \par
$(3)$ The bitension field $\tau_2(\pi)$ satisfies that 
\begin{equation}
\tau_2(\pi):=\int_M\vert\tau(\pi)\vert^2\,v_h=d\,{\mbox{\rm Vol}}(M,h).
\end{equation} 
\end{prop}
\vskip1.2cm\par
By Proposition 3.1, Theorem 1.1 can be proved as follows. 
Assume that $\pi:\,(E,g)\rightarrow (M,h)$ is biharmonic. 
Due to (3.1) in Proposition 3.1, 
we have 
\begin{equation}
\mbox{\rm div}(\tau(\pi)):=\sum_{i=1}^n
(\overline{\nabla}_{e_i'}\tau(\pi))(e_i')=0, 
\end{equation}
where $\{e_i'\}_{i=1}^n$ is a locally defined 
orthonormal frame field on $(M,h)$ 
and we put $n=\dim_{\mathbb R}M$.   
Then, for every $f\in C^{\infty}(M)$, 
it holds that, due to Proposition (3.29) in \cite{U1}, p. 60, for example, 
\begin{equation}
0=\int_Mf\,\mbox{\rm div}(\tau(\pi))\,v_h=
-\int_Mh(\nabla f,\tau(\pi))\,v_h.
\end{equation}
Therefore, we obtain $\tau(\pi)\equiv 0$. \qed
\vskip0.6cm\par
We will prove Proposition 3.1, later. 
Here, we give examples of the line bundles over 
some compact homogeneous K\"{a}hler Einstein manifolds
$(M,h)$: 
\par
\begin{ex}
A generalized flag manifold $G/H$ admits 
a unique K\"{a}hler Einstein metric $h$ (\cite{BH} and \cite{CS}). 
Here, $G$ is a compact semi-simple Lie group, 
and $H$ is the centralizer of a torus $S$ in $G$, i.e., 
$G^{\mathbb C}$ is the complexification of $G$, and $B$ is its 
Borel subgroup. Then, 
$$M=G/H=G^{\mathbb C}/B.$$ 
The Borel subgroup $B$ is written as $B=TN$, where 
 $T$ is a maximal torus of $B$ and 
 $N$ is a nilpotent Lie subgroup of $B$. 
 Every character $\xi_{\lambda}$ of a Borel subgroup $B$ is given as 
 a homomorphism $\xi_{\lambda} :\,B\rightarrow {\mathbb C}^{\ast}={\mathbb C}-\{0\}$ which 
 is written as 
 \begin{align}
 \xi_{\lambda}(tn)=\xi_{\lambda}(t)\qquad (t\in T,\,\,n\in N). 
 \end{align}
 Here $\xi_{\lambda}:\,T\rightarrow U(1)$ is a character of $T$ which 
 is written as 
\begin{equation}
 \xi_{\lambda}(\exp ({\theta_1H_1+\cdots+\theta_{\ell}H_{\ell}}))
 =e^{2\pi\sqrt{-1}(k_1\theta_1+\cdots+k_{\ell}\theta_{\ell})}, \qquad (\theta_1,\,\,\ldots,\,\,\theta_{\ell}\in {\mathbb R}),  
\end{equation}
where $k_1,\,\ldots,\,k_{\ell}$ are non-negative integers, 
and $\ell=\dim T$.  
\par
Note that every character $\xi_{\lambda}$ of a nilpotent Lie group $N$ must be 
 $\xi_{\lambda}(n)=1$ because 
 $\xi_{\lambda}(n)=\xi_{\lambda}(\exp X)=e^{\xi_{\lambda'}(X)}$ where $n=\exp X$ $(X\in {\frak n})$, 
 and $\lambda':\,{\frak t}\rightarrow {\mathbb C}$ is a homomorphism, i.e., $\xi_{\lambda'}(X+Y)=\xi_{\lambda'}(X)+\xi_{\lambda'}(Y)$, $(X,\,Y\in {\frak t})$. 
 Then, there exists $k\in {\mathbb N}$ which 
 satisfies that $\exp (k\, X)=n^k=e$. Then, 
 $e^{k\,\xi_{\lambda'}(X)}=\xi_{\lambda}(n^k)=\xi_{\lambda}(e)=1$. 
 Thus, for every $a\in {\mathbb R}$, 
 $$
 e^{a\,k\,\xi_{\lambda'}(X)}=(e^{k\,\xi_{\lambda'}(X)})^a=1. 
 $$
 This implies that $k\,\xi_{\lambda'}(X)=0$. 
 Thus, $\xi_{\lambda'}(X)=0$ for all 
 $X\in \frak n$, i.e., $\xi_{\lambda'}\equiv 0$. Therefore, we have that $\xi_{\lambda}(n)=e$ 
 ($n\in N$).  We have (3.5). 
 \par
 For every $\xi_{\lambda}$ given by (3.5) and (3.6), 
 we obtain the associated holomorphic vector bundle 
 $E_{\xi_{\lambda}}$ over $G^{\mathbb C}/B$ as 
 $E_{\xi_{\lambda}}:=\{[x,v]\vert (x,v)\in G^{\mathbb C}\times {\mathbb C}\}$, 
 where the equivalence relation $(x,v)\sim (x',v')$ is 
 $(x,v)= (x',v')$ if and only if 
 there exists $b\in B$ such that $(x',v')=(xb^{-1},\xi_{\lambda}(b)v)$,  
 denoted by $[x,v]$, the equivalence class including 
 $(x,v)\in G^{\mathbb C}\times {\mathbb C}$
 (for example, \cite{B}, \cite{TW}). 
\end{ex}
\section{Proof of Proposition 3.1.}
For an Hermitian vector bundle 
$\pi:\,(E,g)\rightarrow (M,g)$ with $\dim_{\mathbb R} E=m$, and 
$\dim _{\mathbb R}M=n$, 
let us recall 
the definitions of the tension field $\tau(\pi)$ and the bitension field 
$\tau_2(\pi)$: 
\begin{equation}
\left\{
\begin{aligned}
\tau(\pi)&=\sum_{j=1}^m\left\{
\overline{\nabla}^h_{e_j}\pi_{\ast}{e_j}
-\pi_{\ast}
\left(
\nabla^g_{e_j}e_j
\right)
\right\},\\
\tau_2(\pi)&=\overline{\Delta}\tau(\pi)
-\sum_{j=1}^mR^h(\tau(\pi),\pi_{\ast}e_j)\pi_{\ast}e_j.
\end{aligned}
\right.
\end{equation}
Then, we have 
\begin{align}
\tau_2(\pi)&:=\overline{\Delta}\tau(\pi)
-\sum_{j=1}^mR^h(\tau(\pi),\pi_{\ast}e_j)\pi_{\ast}e_j\nonumber\\
&=\overline{\Delta}\tau(\pi)-\sum_{j=1}^nR^h(\tau(\pi),e'_j)e'_j\\
&=\overline{\Delta}\tau(\pi)-\mbox{\rm Ric}^h(\tau(\pi)).
\end{align}
Here, recall that $\pi:\,(E,g)\rightarrow (M,h)$ is the Riemannian submersion 
and $\{e_i\}_{i=1}^m$ and $\{e'_j\}_{j=1}^n$ are locally defined orthonormal frame fields on $(E,g)$ and $(M,h)$, respectively, 
satisfying 
that $\pi_{\ast}e_j=e'_j$ $(j=1,\cdots,n)$ and $\pi_{\ast}(e_j)=0$ $(j=n+1,\cdots,m)$. Therefore, we have $(4.2)$ and $(4.3)$ by means of the definition 
of the Ricci tensor field ${\rm Ric}^h$ of $(M,h)$. 
\medskip\par
Assume that $(M,h)$ is a real $n$-dimensional compact K\"{a}hler Einstein manifold 
with ${\mbox{\rm  Ric}}^h=c\,\mbox{\rm Id}$, where $n$ is even.  Then, due to (4.3), 
we have that 
$\pi:\,(E,g)\rightarrow (M,h)$ is biharmonic if and only if 
\begin{equation}
\overline{\Delta}\tau(\pi)=c\,\tau(\pi).
\end{equation}
Since $\langle\tau(\pi),\tau(\pi)\rangle$ is a $C^{\infty}$ function on 
a Riemannian manifold $(M,h)$, we have, 
for each $j=1,\cdots,n$, 
\begin{align}
e'_j\langle \tau(\pi),\tau(\pi)\rangle&=\langle\overline{\nabla}_{e'_j}\tau(\pi),\tau(\pi)\rangle+\langle\tau(\pi),\overline{\nabla}_{e'_j}\tau(\pi)\rangle\nonumber\\
&=2\langle\overline{\nabla}_{e'_j}\tau(\pi),\tau(\pi)\rangle, \\
{e'_j}^{2}\langle\tau(\pi),\tau(\pi)\rangle&=2{e'_j}\langle\overline{\nabla}_{e'_j}\tau(\pi),\tau(\pi)\rangle\nonumber\\
&=2\langle\overline{\nabla}_{e'_j}(\overline{\nabla}_{e'_j}\tau(\pi)),\tau(\pi)\rangle+
2\langle\overline{\nabla}_{e'_j}\tau(\pi), \overline{\nabla}_{e'_j}\tau(\pi)\rangle,\\
\nabla_{e'_j}{e'_j}\langle\tau(\pi),\tau(\pi)\rangle&=2\langle\overline{\nabla}_{\nabla_{e'_j}e'_j}\tau(\pi),\tau(\pi)\rangle. 
\end{align}
Therefore, the Laplacian 
$\Delta_h=-\sum_{j=1}^n(e'_j{}^2-\nabla_{{}_{e'_j}}e'_j)$ acting on 
$C^{\infty}(M)$, so that 
\begin{align}
\Delta_h\left\langle\tau(\pi),\tau(\pi)\right\rangle&=2\sum_{j=1}^n\left\{
-\langle\overline{\nabla}_{e_j'}(\overline{\nabla}_{e'_j}\tau(\pi)),\tau(\pi)\rangle-\langle\overline{\nabla}_{e'_j}\tau(\pi),\nabla_{e'_j}\tau(\pi)\rangle+\langle\overline{\nabla}_{\nabla_{e'_j}}\tau(\pi),\tau(\pi)\rangle
\right\}\nonumber\\
&=2\big\langle-\sum_{j=1}^n\big\{\overline{\nabla}_{e'_j}\overline{\nabla}_{e'_j}-\overline{\nabla}_{\nabla_{e'_j}e'_j}
\big\}\tau(\pi),\tau(\pi)\big\rangle
-2\sum_{j=1}^n\big\langle
\overline{\nabla}_{e'_j}\tau(\pi),\overline{\nabla}_{e'_j}\tau(\pi)\big\rangle\nonumber\\
&=2\big\langle\overline{\Delta}\tau(\pi),\tau(\pi)\big\rangle-2\sum_{j=1}^n\langle\overline{\nabla}_{e'_j}\tau(\pi),\overline{\nabla}_{e'_j}\tau(\pi)\rangle\\
&\leq 2\big\langle\overline{\Delta}\tau(\pi),\tau(\pi)\big\rangle, 
\end{align}
because of  
$\langle\overline{\nabla}_{e'_j}\tau(\pi),\overline{\nabla}_{e'_j}\tau(\pi)\rangle\geq0$, $(j=1,\cdots,n)$. 
\par
If $\pi:\,(E,g)\rightarrow (M,h)$ is biharmonic, 
due to (4.4), $\overline{\Delta}\tau(\pi)=c\,\tau(\pi)$, 
the right hand side of (4.8) coincides with 
\begin{align}
(4.8)&=
2c\,\langle\tau(\pi),\tau(\pi)\rangle-
2\sum_{j=1}^n\langle\overline{\nabla}_{e'_j}\tau(\pi),\overline{\nabla}_{e'_j}\tau(\pi)\rangle\\
&\leq2c\,\langle\tau(\pi),\tau(\pi)\rangle.
\end{align} 
Remember that 
due to M. Obata's theorem, (see Proposition 4.1 below), 
\begin{equation}
\lambda_1(M,h)\geq 2c, 
\end{equation}
since 
$\mbox {\rm Ric}_h=c\,\mbox{\rm Id}$, 
 And the equation in (4.11) holds, i.e., 
$\lambda_1(M,h)=2c$ and 
 \begin{equation}
 \Delta_h\,\langle\tau(\pi),\tau(\pi)\rangle=
 2c\,\langle\tau(\pi),\tau(\pi)\rangle
 \end{equation}
 holds. Then, (4.12) implies that
 the equality in the inequality (4.11) holds. We have that 
 \begin{equation}
\sum_{j=1}^n\langle\overline{\nabla}_{e'_j}\tau(\pi), 
\overline{\nabla}_{e'_j}\tau(\pi)\rangle=0, 
\end{equation}
which is equivalent to that
\begin{equation}
\overline{\nabla}_X\tau(\pi)=0\qquad\qquad (\forall X\in {\frak X}(M)).
\end{equation}
Due to (4.15), for every $X\in {\frak X}(M)$, 
\begin{equation}
X\,\langle\tau(\pi),\tau(\pi)\rangle=2\,\langle\overline{\nabla}_X\tau(\pi),\tau(\pi)\rangle=0.
\end{equation}
Therefore, the function 
$\langle\tau(\pi),\tau(\pi)\rangle$ on $M$ is a constant function on $M$. 
Therefore, it implies that 
the right hand side of (4.12) must vanish. Thus, $c=0$ or 
$\tau(\pi)\equiv 0$.  If we assume that $\tau(\pi)\not\equiv 0$, then 
by (4.12), it must hold that $2c=0$.  
Then, $\overline{\Delta}\tau(\pi)=c\,\tau(\pi)=0$, so that 
$\tau(\pi)\equiv 0$ due to (4.4). 
\vskip0.6cm\par
Let $\lambda_1(M,g)$ be the first eigenvalue of the Laplacian 
$\Delta$ of a compact Riemannian manifold $(M,g)$. 
Recall the theorem of M. Obata: 
\begin{prop}(cf. \cite{U1}, pp. 180, 181 )
Assume that $(M,g)$ is a compact K\"{a}hler manifold, 
and the Ricci transform $\rho$ of $(M,g)$ satisfies that 
\begin{equation}
g(\rho(u),u)\geq \alpha\,g(u,u),\qquad(\forall u\in T_xM), 
\end{equation} 
for some positive constant $\alpha>0$. 
Then, it holds that 
\begin{equation}
\lambda_1(M,g)\geq 2\,\alpha.
\end{equation}
If the equality holds, then $M$ admits a non-zero holomorphic 
vector field. 
\end{prop}
Thus, we obtain Proposition 3.1, and the following theorem 
(cf. Theorem 1.1): 
\begin{thm}
Let $\pi:\,(E,g)\rightarrow (M,g)$ be an Hermitian vector over 
a compact K\"{a}hler Einstein manifold $(M,h)$. 
If $\pi$ is biharmonic, then it is harmonic.  
\end{thm}
\newpage
%

\end{document}